\documentclass[12pt]{article}
\usepackage{amssymb}
\usepackage{graphicx}
\usepackage{latexsym}
\def\part#1{\frac{\partial\phantom{q}}{\partial#1}}

\newenvironment{rmks}{\begin{trivlist}\item[]{\bf Remarks:} }
{\end{trivlist}}

\newenvironment{prf}{\begin{trivlist}\item[]{\bf Proof:} }
{\hfill $\Box$ \end{trivlist}}

\newtheorem{thm}{Theorem}

\newtheorem{prp}[thm]{Proposition}

\newcommand{\lie}[1]{\mathfrak{#1}}
\def\End{\mathop{\rm End}\nolimits}
\def\Hom{\mathop{\rm Hom}\nolimits}

\def\tr{\mathop{\rm tr}\nolimits}

\def\ad{\mathop{\rm ad}\nolimits}

\newcommand{\C}{\mathbf{C}}

\newcommand{\Z}{\mathbf{Z}}

\newcommand{\PP}{{\rm P}}

\newcommand{\OO}{{\mathcal O}}

\begin{document}
\title{Remarks on the intersection of two quadrics}
\author{Nigel Hitchin\\[5pt]}
\maketitle

\centerline{{\it Dedicated to G\'erard Laumon}} 
\section{Introduction}
Two papers, one by G\'erard Laumon \cite{L} and one by the author \cite{H1} appeared independently in a collected volume of the Duke Mathematical Journal dedicated to Yuri Manin. Both concern the cotangent bundle of a space of holomorphic principal  $G$-bundles $P$ over a curve $C$, where a cotangent vector is interpreted as a section $\Phi$ of the vector bundle $\ad(P)\otimes K$. Laumon in his paper showed that the subvariety for which $\Phi$ is nilpotent is Lagrangian, and the author in his introduced a space of Poisson-commuting functions yielding an integrable system,   and identified the generic fibre (for the classical groups)  with an open set in an abelian variety. This is also Lagrangian and in this context Laumon's  nilpotent cone is the special fibre over the origin.  

The introduction of a  concrete example eluded the author in \cite{H1}, and one purpose of this article is to address that omission. The simplest situation consists of stable rank 2 vector bundles $E$, where $\Lambda^2E$ is fixed of odd degree, over a curve $C$ of genus $2$. The  moduli space   was identified by Newstead \cite{New} as the 3-dimensional intersection of two quadrics.  This case has, however, been generalized in a recent paper by Beauville at al. \cite{Betal}, producing an explicit formula for  an  integrable system on the cotangent bundle of an intersection of quadrics of arbitrary dimension $n$. In this article we  take that formula, reinterpret it, and discuss some of the issues, notably the concept of {\it very stable points} introduced in Laumon's article, in this particular  context. 

The very stable bundles are those for which there is no non-zero nilpotent $\Phi$ (which we shall now call a Higgs field). It is an open set in the moduli space of stable bundles and its complement, the {\it wobbly} locus, has been the subject of several recent investigations because of its role in aspects of the Geometric Langlands programme, especially in the approach of Donagi et al. \cite{Don}. It is a concept which makes sense for any integrable system defined on a cotangent bundle. We show here that, for the intersection $X=Q\cap Q_1$ of dimension $n$, 
the quotient  by $\Z_2^{n+2}$ is isomorphic to $\PP^n$ realized as the $n$-fold symmetric product of $\PP^1$ and the very stable points map to the $n$-tuples of {\it distinct} points. 
We then pursue the analogy with moduli of bundles a little further and observe that the formula in \cite{Betal} for the Poisson-commuting functions directly defines a family of commuting differential operators on the square root  of the canonical bundle of $X$ and hence offers the opportunity to explore a version of the analytic Langlands correspondence as in \cite{EFK}.

All of the foregoing capitalizes on the explicit form of the functions generating the integrable system by expressing them in terms of quasi-parabolic rank $2$ vector bundles on $\PP^1$. However a very recent paper \cite{BHL} shows that, in odd dimensions, this is the usual integrable system defined on the cotangent bundle of the moduli space of semi-stable twisted $Spin(2g)$ bundles over a hyperelliptic curve $C$ of genus $g$, invariant under the hyperelliptic involution $\tau$, following the description of Ramanan \cite{Ram}. A related description in the even case is also given. The remarkable feature in this interpretation is that the Higgs field $\Phi$, locally taking values in the Lie algebra $\lie{so}(2g)$,  in fact has rank $2$. In the final section we reveal the link with our description above. This interpretation now offers the opportunity to produce a moduli space of stable Higgs bundles in which Laumon's nilpotent cone compactifies.

\section{A symplectic quotient}\label{symp}
Let $V\cong \C^2$ be a 2-dimensional symplectic vector space with skew form $\langle v,w\rangle$ and consider the  action of $SL(2,\C)$ on $V\otimes \C^{N+1}$. The moment map is 
$$m(v_1,\dots, v_N, v_0)= \sum_{i=0}^N v_i\otimes v_i$$
where we identify the lie algebra of $SL(2,\C)$ with $S^2V$. Each $\langle v_i,v_j\rangle$ is an invariant function and hence defined on the symplectic quotient $m^{-1}(0)/SL(2,\C)$. Since 
$$v_0\otimes v_0=-\sum_{i=1}^N v_i\otimes v_i$$
it follows that $ \langle v_i,v_0\rangle ^2+\cdots +  \langle v_N,v_0\rangle ^2=0$. Put $x_i= \langle v_i,v_0\rangle$. 

Take a symplectic basis $e_1,e_2$ of $V$ with $v_0=e_2$, then $v_i=x_ie_1+y_ie_2$ and setting $m=0$ gives 
\begin{equation}\label{xy}
\sum_{i=1}^N x_i^2=0,\quad \sum_{i=1}^Ny_i^2=-1,\quad  \sum_{i=1}^N x_iy_i=0
\end{equation}
 and the   stabilizer of $e_2$ acts as $y_i\mapsto y_i+tx_i$. This leaves $x\wedge y\in \Lambda^2\C^N$ unchanged where $x=(x_1,\dots, x_N), y=(y_1,\dots, y_N)$  and for $x\wedge y\ne 0$ we obtain a symplectic manifold which is a coadjoint orbit of $SO(N)$ acting on  $\C^N$. It consists of the non-null cotangent vectors of the quadric $Q: \{x_1^2+\cdots +x_N^2=0\}\subset \PP^{N-1}$. 
 
 The function $\langle v_i,v_j\rangle=x_iy_j-x_jy_i$ is the moment map for the canonical lift of the vector field $X_{ij}=x_i\partial_j-x_j\partial_i$ on $Q$ to $T^*Q$.
 \section{Parabolic bundles}
 Consider now the meromorphic Higgs field on the trivial bundle over $\PP^1$ defined by 
\begin{equation}\label{phi}
\Phi=\sum_{i=1}^N\frac{v_i\otimes v_i}{z-\mu_i}dz. 
\end{equation}
 This has nilpotent residue $v_i\otimes v_i$ at $z=\mu_i$ and since $v_1\otimes v_1+\dots +v_0\otimes v_0=0$ and the sum of residues is zero, we have a  nilpotent residue $v_0\otimes v_0$ also at infinity. 
 
 The parabolic version of the integrable system \cite{H1} has Poisson commuting functions given by the coefficients of 
 $$\tr \Phi^2=-\sum_{i,j}\frac{\langle v_i,v_j\rangle^2}{(z-\mu_i)(z-\mu_j)}dz^2$$
 or equivalently, by taking the residue at $z=\mu_i$ the functions 
 $$f_i=\sum_{j\ne i}\frac{(x_iy_j-x_jy_i)^2}{(\mu_i-\mu_j)}.$$
 Written as a symmetric tensor this is 
\begin{equation}\label{si}
s_i=\sum_{j\ne i}\frac{(x_i\partial_j-x_j\partial_i)^2}{(\mu_i-\mu_j)}.
 \end{equation}
 Now consider $q_1= \mu_1x_1^2+\cdots +\mu_N x_N^2$ and take the inner product of the symmetric tensor $s_i$ in (\ref{si}) with $dq_1$. We obtain
 $$\sum_{j\ne i}4x_ix_j(\mu_j-\mu_i)\frac{(x_i\partial_j-x_j\partial_i)}{(\mu_i-\mu_j)}=-4x_i^2\sum_{j\ne i}x_j\partial_j+4x_i\partial_i\sum_{j\ne i}x_j^2=-4x_i^2\sum_{j=1}^Nx_j\partial_j$$
 using $x_1^2+\cdots +x_N^2=0$. But this is a multiple of the Euler vector field on $\C^N$ and hence is zero in $\PP^{N-1}$. We conclude that $s_i$  defines a symmetric tensor on the intersection $X$ of $Q$ with $Q_1$ defined by $q_1=0$ and (\ref{si})  is the formula in \cite{Betal} (Proposition 7.4). 
 
Returning to the Higgs field $\Phi$,  observe that as $z\rightarrow \infty$,
\begin{equation}\label{expand}
\Phi=- v_0\otimes v_0\frac{dz}{z}+\sum_{i=1}^N \mu_i v_i\otimes v_i \frac{dz}{z^2}+\cdots =\phi_0 \frac{dz}{z}+ \phi_1\frac{dz}{z^2}+\cdots 
\end{equation}
 and the equation $q_1=0$ is equivalent to $\tr \phi_0\phi_1=0$, which implies that $v_0$ in the kernel of $\phi_0$ is  also an eigenvector of $\phi_1$. We may therefore apply a Hecke transformation at infinity using the distinguished subspace defined by $v_0$  to remove the singularity of $\Phi$ at the expense of considering $E=\OO\oplus \OO(-1)$ with $v_0$ now interpreted as a section of the unique trivial subbundle of $E$.
 
 \begin{rmks} 
 
 \noindent 1. The spectral curve $S$ for  a generic $\Phi$ is the hyperelliptic curve $y^2+\det \Phi=0$, and a line bundle on $S$ produces by direct image the bundle $E$ and Higgs field $\Phi$. Since $\Phi$ only determines the $v_i$ up to sign, the fibre in $T^*X$ of the integrable system is an unramified covering of an open set of the Jacobian of $S$. It is described in \cite{Betal} as a \emph {quotient} of the Jacobian but the two are related by the map of divisor classes $x\mapsto 2x$. 
 
 \noindent 2. 
 Using the basis $e_1, e_2$ the Higgs field $\Phi$ is given by 
$$\Phi(e_1) = -\sum_{i=1}^N\frac{x_iy_i}{z-\mu_i} e_1- \sum_{i=1}^N\frac{y^2_i}{z-\mu_i} e_2$$
$$\Phi(e_2) = \sum_{i=1}^N\frac{x_i^2}{z-\mu_i} e_1+ \sum_{i=1}^N\frac{x_iy_i}{z-\mu_i} e_2$$
 and in this form is recognizable as the Garnier system or the classical Gaudin system. The reader may see this in a more  general context in the survey lectures \cite{EL}.
 
 \end{rmks}

 \section{Polynomials}
 We now view   the Higgs field as $\Phi:E\rightarrow E\otimes K(D)$  where $E\cong \OO\oplus \OO(-1)$ and $D$ is the divisor of the points $z=\mu_i$. The trivial subbundle is unique and hence $\Phi$ determines another divisor  consisting of the points $a_i\in \PP^1$ at which $\Phi$ preserves $\OO\subset E$. 
    In the formulas above this means $\Phi(v_0)=\lambda v_0$ or $\langle \Phi (v_0),v_0\rangle=0$ or 
     \begin{equation}\label{aux}
     0=\langle \Phi (v_0),v_0\rangle=\sum_{i=1}^N\frac{\langle v_i,v_0\rangle^2}{z-\mu_i}=\sum_{i=1}^N\frac{x_i^2}{z-\mu_i}.
     \end{equation}    
            From the expansion at infinity (\ref{expand})   we see that the leading terms    vanish when $x_1^2+\cdots +x_N^2=0= \mu_1x_1^2+\cdots + \mu_Nx_N^2$ and so clearing the denominators this gives a polynomial $p(z)$ in $z$ of degree $n=N-3$ with roots $z=a_1,\dots, a_n$.

            We calculate the eigenvalue $\lambda_k$ at $z=a_k$ directly, supposing $a_k\ne \mu_j$ for any $j$:
            $$\sum_{i=1}^N\frac{\langle v_i,v_0\rangle v_i}{a_k-\mu_i}=\lambda_k v_0$$
            and hence
            $$\sum_{i=1}^N\frac{\langle v_i,v_0\rangle \langle v_i, v_j\rangle }{a_k-\mu_i}=\lambda_k \langle v_0,v_j\rangle $$
               for all $j$, or 
               $$ \sum_{i=1}^N\frac{x_i (x_iy_j-x_jy_i)}{a_k-\mu_i}=-\lambda_k x_j $$
               which, using (\ref{aux}) gives 
              \begin{equation}\label{eigen}
              \lambda_k=  \sum_{i=1}^N\frac{x_iy_i}{a_k-\mu_i}
              \end{equation}
The Poisson-commuting functions are, in this parabolic context, the coefficients of $\tr \Phi^2\in H^0(\PP^1,K^2(2D))=H^0(\PP^1,\OO(2N-4))$. Since $\Phi$ is nilpotent at $z=\mu_i$ this reduces to $H^0(\PP^1,\OO(N-4))$ of dimension $n=N-3=\dim X$. The integrable system is a map $f:T^*X\rightarrow H^0(\PP^1,\OO(n-1))$.

Evaluation of a polynomial of degree $(n-1)$ at $n$ distinct points is an alternative basis to using the coefficients of powers of $z$, so if we assume the  $a_i$ are distinct, we may use these points to describe the quadratic functions on the cotangent bundle. But at $z=a_k$, $\Phi (v_0)=\lambda_k v_0$ and $\tr \Phi^2=-\lambda_k^2$.  It follows from (\ref{eigen}) that we obtain the $n$ Poisson-commuting functions 
\begin{equation}\label{linear}
f_k=\left(\sum_{i=1}^N\frac{x_iy_i}{a_k-\mu_i}\right)^2 = \ell_i(y)^2
\end{equation}
We have  $n$ linear forms $\ell_i(y)$ where also $\sum_i x_iy_i=0=\sum_i \mu_i x_iy_i$.  Suppose $\ell_i(y)=0$ for all $i$, then the  determinant of the following matrix must vanish:
\begin{equation}\label{matrix}
\pmatrix{x_1(a_1-\mu_1)^{-1} & x_2(a_1-\mu_2)^{-1} & \cdots & x_N(a_1-\mu_N)^{-1}\cr
x_1(a_2-\mu_1)^{-1} & x_2(a_2-\mu_2)^{-1} & \cdots & x_N(a_2-\mu_N)^{-1}\cr
\cdots & \cdots & \cdots &\cr
x_1 & x_2 & \cdots & x_N\cr
\mu_1 x_1 & \mu_2 x_2 & \cdots & \mu_N x_N}
\end{equation}
and this evaluates to 
$$\pm \prod_{i=1}^N  \frac{x_i}{p(\mu_i)}\prod_{j<k}(a_j-a_k)\prod_{\ell<m}(\mu_{\ell}-\mu_n)$$
where the $a_j$ are the roots of $p(z)=0$ where $p(z)=x_1^2(z-\mu_2)\dots (z-\mu_N)+\cdots$ so $p(\mu_1)=x_1^2(\mu_1-\mu_2)(\mu_1-\mu_3)\dots (\mu_1-\mu_N)$ etc. 

Hence the determinant is non-zero since the $\mu_i$  are distinct for a smooth intersection and the  $a_j$ distinct by assumption. This means that the $\ell_i$ are linearly independent and $\ell_i^2$ span the $n$-dimensional  space of functions for the integrable system. 
\begin{rmks} 

\noindent 1. In the context of mirror symmetry for Higgs bundles the notion of multiplicity algebra was introduced in \cite{H}. For a  principal  $G$-bundle this is the algebra with relations defined by the invariant polynomials of $G$ on the cotangent space.  Even in rank $2$, these can be quite complicated (see \cite{H3} for example) but in the above case we have seen that it is a sum of squares of linear functions if the $a_i$ are distinct and disjoint from the $\mu_j$. 

\noindent 2. From the point of view of the spectral curve, the vector bundle  $E$ is the direct image of a line bundle $L$ of degree $n$ under the projection $\pi:S\rightarrow \PP^1$. Since the vector bundle is $\OO\oplus \OO(-1)$ there is a canonical section $s$ of $L$ 
corresponding to $v_0$ spanning the trivial subbundle. The image in $\PP^1$ of the divisor of $s$  is $a_1+\dots +a_n$.

 \noindent 3. In the integrable systems community, passing to the coordinates $a_1,\dots, a_n$ is known as separation of variables as in \cite{Skl}, or the survey article \cite{Hurt}. In our case the vector bundle $E$ has a canonical section yielding these points, rather than the choice involved in the general case. 

 \noindent 4. In \cite{At}, Atiyah described projective bundles over a curve $C$ of genus $2$ in terms of vector bundles $\OO\rightarrow E\rightarrow L$ where $L$ has degree $1$. There is a unique section of $LK$ which lifts to $E\otimes K$ and its divisor consists of three points on $C$ which project to our points $a_1,a_2,a_3$ in $\PP^1$ under the quotient map of the hyperelliptic involution.

\end{rmks}
\section{Very stable points} 
Given an integrable system on the cotangent bundle of a manifold $M$ one may say that a point in $M$ is very stable if there are no cotangent vectors for which all the functions of the integrable system vanish. In the case of a principal $G$-bundle this means a Higgs field $\Phi$ such that all invariant polynomials on $\lie{g}$ vanish, and hence $\Phi$ is nilpotent everywhere. For the intersection of quadrics, what we showed in the previous section was that if all the functions vanish, and $x_i\ne 0$ for any $i$ then  $\ell_i(y)=0$ for all  $y$. This implies $y=0$ if the $a_j$ are distinct, and then $(x_1,\dots, x_N)$ is a very stable point. 

As a consequence, if we have a nilpotent Higgs field, then there must be a multiple zero of $p(z)$, and we can see this in general  directly: since $E=\OO\oplus \OO(-1)$ we write 
\begin{equation}\label{phimatrix}
\Phi=\pmatrix{ b & a\cr
                         c & -b}
                         \end{equation}
                         where $c\in H^0(\PP^1,\OO(n)), b\in H^0(\PP^1,\OO(n+1)), a\in H^0(\PP^1,\OO(n+2))$, so that $c$ vanishes when $\OO$ is preserved, hence $c$ is essentially $p(z)$. Then $\Phi$ is nilpotent if $b^2+ac=0$ and at a zero $a_i$ of $c$, $b=0$ so if $a\ne 0$ then $c=-a^{-1}b^2$ has a double zero. The full result is the following:

                         \begin{prp} A point $(x_1,\dots, x_N)$ on the intersection of quadrics $Q\cap Q_1$ is very stable with respect to the integrable system if and  only if the polynomial
                         $$p(z)=\sum_{i=1}^N x_i^2 (z-\mu_1)(z-\mu_2)\dots \widehat{(z-\mu_i)}\dots (z-\mu_N)$$
                         has distinct roots (including $z=\infty$).
                          \end{prp} 
                          \begin{prf} Framed in terms of Higgs fields acting on $\OO\oplus \OO(-1)$ we can transform by the   action of $PGL(2,\C)$ on $\PP^1$.  Therefore $z=\infty$ has no distinguished role, and we may assume that the $a_i$, roots of $p(z)=0$ are finite. In the previous section we assumed $a_i\ne \mu_j$ so we should consider the case $a_1=\mu_1$ for example, then $x_1^2(\mu_1-\mu_2)\dots (\mu_1-\mu_N)=0$ and so $x_1=0$ and we are in the situation of an intersection of quadrics of dimension $n-1$. Induction on $n$ then incorporates this case. Appealing to the previous section we see that if the $a_i$ are distinct  we have a very stable point. We now need the converse. 
                          
                          Suppose now $a_2=a_1$ and $z=a_1$ is a double zero of $p(z)$, with $a_1, a_3,\dots a_n$  distinct. Then the matrix (\ref{matrix}) is singular and we have a nonzero $y$ in the cotangent space such that $\ell_i(y)=0$ for $i\ne 2$.               Now choose a basis for $H^0(\PP^1,\OO(n-1))$ by evaluation at $a_1,a_3,\dots a_n$ and evaluation of the derivative at $a_1$. In the matrix form (\ref{phimatrix}) near $z=a_1$ we have $c=(z-a_1)^2(c_0+\dots), b=(z-a_1)(b_0+\dots)$ so that $\tr\Phi^2=(z-a_1)^2(d_0+\dots)$ which vanishes at $z=a_1$ together with its first derivative. Hence all functions of the integrable system vanish and the Higgs field is nilpotent. The argument can be modified for several multiple zeros.
                             \end{prf}
                          \begin{rmks} 
                          
                          \noindent 1. As shown in  \cite{PP}, a vector bundle $E$ on a curve is very stable if and only if the map from the cotangent space $H^0(C, \End E\otimes K)$   to the base $\C^n$ of the integrable system is proper.  In our case, the map is a linear isomorphism of $\C^n$ followed by the map $(y_1,\dots, y_n)\mapsto (y_1^2,\dots, y_n^2)$ which is clearly proper. 
                          
                          \noindent 2. The complement of the very stable points is the inverse image of the discriminant hypersurface given by the resultant of $p$ and $p'$ under the map $X\rightarrow \PP^n$ defined by $x_i\mapsto x_i^2$ $1\le i\le N$. In the case $n=3$, concerning stable bundles on a genus 2 curve, this is observed in \cite{Don}.

                          \end{rmks} 
\section{Differential operators}
     The authors of \cite{Betal} present the Poisson-commuting functions in the form  $$s_i=\sum_{j\ne i}\frac{(x_i\partial_j-x_j\partial_i)^2}{\mu_i-\mu_j}$$     considered as sections of the symmetric power $S^2T$ of the tangent bundle of $X$, but this can also be interpreted as  a second-order differential operator.     In fact $X_{ij}=x_i\partial_j-x_j\partial_i$ is the vector field on the quadric $Q$ defining rotation in the $x_i,x_j$-plane, these elements providing the standard basis for the lie algebra $\lie{so}(N)$. From this viewpoint $s_i$ is a differential operator $\Delta_i$ acting on local sections of any homogeneous vector bundle over $Q$. Take the line bundle $\OO(k)$ from the embedding $Q\subset \PP^{N-1}$: we would like  to define $\Delta_i$ as an operator on local sections on $Q\cap Q_1$, since we have already seen that its symbol is defined on the intersection.
     
     So consider its action on $fq_1$ where $f$ is a local section of $\OO(k-2)$, which we consider as a function $f(x)$ homogenous of degree $(k-2)$.  We have 
 $$\Delta_i(fq_1)=\Delta_i(f)q_1+2\sum_{j\ne i}\frac{(x_i\partial _jf-x_j\partial_if)(x_i\partial _jq_1-x_j\partial_iq_1)}{\mu_i-\mu_j}+f\Delta_iq_1.$$
 Now $(x_i\partial _jq_1-x_j\partial_iq_1)=2(\mu_j-\mu_i)x_ix_j$ and so the middle term may be written as $$-4\sum_{j\ne i}x_ix_j(x_i\partial _jf-x_j\partial_if)=-4x_i^2\sum_{j\ne i}(x_j\partial _jf)+4x_i\partial_if\sum_{j\ne i}x_j^2=-4x_i^2(k-2)f$$
 using the homogeneity of $f$ and $x_1^2+\cdots +x_N^2=0$. 
 And, using $(x_i\partial _jq_1-x_j\partial_iq_1)=2(\mu_j-\mu_i)x_ix_j$ again,  
 $$\Delta_iq_1=-2\sum_{j\ne i}(x_i\partial _j-x_j\partial_i)x_ix_j=-2\sum_{j\ne i}(x_i^2-x_j^2)=-2(N-1)x_i^2 -2x_i^2=-2N x_i^2$$
 Hence 
 $\Delta_i(fq_1)=\Delta_i(f)q_1+x_i^2(-4(k-2)-2N)f$ which means that if we take $k=-(N-4)/2$ then $\Delta_i$ is a well-defined holomorphic operator on $X$, since divisibility by $q_1$  is preserved. 
 
 Now, as an intersection of two quadrics, $K_X\otimes \OO(-2)\otimes \OO(-2)\cong K_{\PP^{N-1}}\cong \OO(-N)$ so $K_{X}\cong \OO(-(N-4))$ and $\OO(k)$ is a square root of the canonical bundle.  This provides a setting analogous to \cite{EFK} where the operators and their conjugates act on global $C^{\infty}$ sections of $K^{1/2}_{X}\otimes \bar K^{1/2}_{X}$, where there is a natural $L^2$ inner product. 
 
 The operators $\Delta_i,\Delta_j$ commute as can be seen by a direct calculation: $\Delta_i$ is of the form 
 $$\sum_{j\ne i}\frac{\Omega_{ij}}{\mu_i-\mu_j}$$
 with $\Omega_{ij}=\Omega_{ji}$ exactly as in the KZ-equation and these operators commute if the Kohno-Drinfeld relations hold: $[\Omega_{ij},\Omega_{k\ell}]=0$ if the indices are distinct and otherwise $[\Omega_{ij}, \Omega_{ik}+\Omega_{jk}]=0$. 
 
 However, with $\Omega_{ij}=X_{ij}^2$, the vector fields $X_{ij}, X_{ik}, X_{jk}$ form a basis for a copy of $\lie{so}(3)\subset \lie{so}(N)$ and then the Casimir $X^2_{ij}+X^2_{ik}+X_{jk}^2=\Omega_{ij}+\Omega_{ik}+\Omega_{jk}$ commutes with everything, in particular $\Omega_{ij}$, giving the second relation. The first is clear since  $X_{ij}, X_{k\ell}$  are rotations in orthogonal planes. 
 
\section{$X$ as a moduli space}
     Ramanan's paper \cite{Ram} identifies the intersection of quadrics $X$ of odd dimension $n=2g-1$ as  the moduli space of semi-stable $Spin(2g)$-bundles on a hyperelliptic curve $C$ of genus $g$, invariant by the hyperelliptic involution $\tau$. The authors of \cite{BHL} show that the original integrable system in \cite{H1}, restricted to the fixed point set of $\tau$,   is in fact equivalent to the one introduced in the earlier paper \cite{Betal}. Here the curve $C$ is the double covering of $\PP^1$ branched over the points $z=\mu_i$. We explain now the link with our rank 2 version, following  the remarkable result from \cite{BHL} that the $Spin(2g)$-Higgs field has rank $2$. 
     
   Higgs fields take values in the adjoint representation, so for their consideration it is enough to look at the associated   orthogonal bundle, rather than the spin bundle. The work of Bhosle \cite{Bhosle} shows that a $\tau$-invariant rank $2g$ orthogonal bundle on $C$ is equivalent to a degenerate orthogonal bundle on $\PP^1$, a bundle $W$ with a homomorphism $W\rightarrow W^*(1)$ which drops rank at the points $z=\mu_i$.
   
    In Ramanan's construction the projective line $\PP^1$ is viewed as the base of the pencil of quadrics $zq-q_1=0$ in $\PP^{N-1}$, so a point in $X$ defines a one-dimensional subspace $L\subset \C^N$, isotropic with respect to the quadratic form  $q_z=zq-q_1$. Then the degenerate orthogonal bundle is given by $W=L^{\perp}/L$ of rank $N-2=2g$ where orthogonality is defined using $q_z$, hence an inner product with values in $\OO(1)$. 
    
    We use the  coordinates $(x_i,y_i)$ as in Section \ref{symp} so that $x=(x_1,\dots, x_N)$ is a point on $X$ generating  the subspace $L\subset \C^N$. Then 
    $$L^{\perp}=\{u\in \C^N: \sum_{i=1}^N(z-\mu_i)x_iu_i=0\}.$$
  This clearly contains the subspace $ \sum_{i=1}^Nx_iu_i=0=\sum_{i=1}^N\mu_ix_iu_i$ defining tangents to $X$, but more invariantly since $T_x\PP^{N-1}=\Hom(L,\C^N/L)$, we have 
  $$L\otimes TX_x\subset L^{\perp}/L=W.$$ 
     Recall now that  $L$ is fixed but $L^{\perp}$ varies over $\PP^1$ with the quadratic form $q_z$, so $L\otimes TX_x$ is a trivial subbundle of $W$.
      Since $q_z$ has rank $N-1$  at each point $z=\mu_i$   we have $(\Lambda^{N-2}W^*)^2(N-2)\cong \OO(N)$ and so $\Lambda^{N-2}W\cong \OO(-1)$. It follows (as in \cite{BHL}) that $W$ is an extension 
      $$0\rightarrow L^*\otimes TX_x\rightarrow W\rightarrow \OO(-1)\rightarrow 0,$$
     where the trivial subbundle is unique, but   the extension splits since $H^1(\PP^1,\OO(1))=0.$
     
     Consider now the matrix 
     $$A_{ij}=\frac{x_iy_j-y_ix_j}{z-\mu_i}.$$
     This is clearly skew-adjoint with respect to the inner product $q_z$. Furthermore, since $\sum_{j=1}^Nx_j^2=0=\sum_{j=1}^Nx_jy_j$, the vector $x$ lies in the kernel, namely $L$,  and so $A$ preserves $L$ and $L^{\perp}$. It thus defines a skew-adjoint meromorphic endomorphism of $W=L^{\perp}/L$. Restrict to $u\in L^*\otimes TX_x$, and we obtain
     $$\sum_{j=1}^N\frac{x_iy_ju_j}{z-\mu_i}-\frac{y_ix_ju_j}{z-\mu_i}=\sum_{j=1}^Ny_ju_j\left(\frac{x_i}{z-\mu_i}\right)$$
   so the intersection $U$ of $L^*\otimes TX_x$ with the kernel is given by $\sum_{j=1}^Nu_iy_i=0$ which is the annihilator in $ L^*\otimes TX_x$ of the cotangent vector defined by $y$. Then each cotangent vector $y$ to $X$ at $x$  defines a meromorphic rank 2 skew-adjoint endomorphism $A$ of $W$ and, appealing to \cite{Bhosle},  $Adz$ defines a $\tau$-invariant Higgs field on $C$.

  \begin{prp} The  Higgs field $$\Psi=Adz: W/U\rightarrow W/U\otimes K$$ is equivalent to $\Phi$ in equation (\ref{phi}).
  \end{prp} 
\begin{prf} 
If $v\in \C^N$ satisfies $\sum_{i=1}^N v_ix_i=0$, then 
$$v_z=\left(\frac{v_1}{z-\mu_1},\dots, \frac{v_N}{z-\mu_N}\right)$$
lies in $L^{\perp}$ and so has an image in $W=L^{\perp}/L$. In particular we have $x_z,y_z\in L^{\perp}$. 
We calculate 
$$A(x_z)=\sum_{i=1}^N\frac{x_iy_i}{z-\mu_i} x_z-\sum_{i=1}^N\frac{x_i^2}{z-\mu_i} y_z$$
$$A(y_z)=\sum_{i=1}^N\frac{y^2_i}{z-\mu_i} x_z-\sum_{i=1}^N\frac{x_iy_i}{z-\mu_i}y_z$$

Comparing this with the explicit expression for $\Phi$ in Section 3  we see that 
 these coincide if $e_1=y_z,e_2=-x_z$. 

On the face of it, $x_z,y_z$ appear to be singular but it is their image in $W/U$ which we need.  The $z^{-1}$ term in the expansion of $\sum_{i=1}^N x_i(v_z)_i$ as $z\rightarrow \infty$ is zero by definition and the next term is $\sum_{i=1}^N\mu_ix_iv_i$ which represents the map $W\rightarrow \OO(-1)$. This vanishes for $v=x$ which means $x_z$ maps to the trivial subbundle $\OO\subset W/U$, and indeed $e_2$ was defined this way. Then $x_z$ together with the image of $y_z$ gives a local basis. 
\end{prf}

  \vskip 2cm
  Mathematical Institute
  
  University of Oxford 
  
  Woodstock Road
  
  Oxford OX2 6GG
  
  UK

\end{document}